\documentclass[12pt]{iopart}
\usepackage{times}
\usepackage{graphicx}

\usepackage{amssymb}
\newcommand\RR{{\mathbb R}}

\newcommand\ZZ{{\mathbb Z}}
\newcommand\TT{{\mathbb T}}
\newtheorem{thm}{Theorem}[section]
\newtheorem{lemma}[thm]{Lemma}
\newtheorem{cor}[thm]{Corollary}
\newtheorem{prop}[thm]{Proposition}

\newtheorem{rem}[thm]{Remark}

\newenvironment{prf} {{\bf Proof.}}{\hfill \framebox(8,8){ } }

\begin{document}

\title{Admissible vectors and traces on the commuting algebra}
\author{Hartmut F\"uhr}
\address{Institute for Biomathematics and Biometry, GSF National Research
 Center for Environment and Health, Ingolst\"adter Landstra\ss{}e 1,
 D-85764 Neuherberg}
\begin{abstract}
Given a representation of a unimodular locally compact group,
we discuss criteria for associated coherent state expansions
in terms of the commuting algebra. It turns out that for those
representations that admit such expansions there exists a unique
finite trace on the commuting algebra such that the admissible
vectors are precisely the tracial vectors for that trace.
This observation is immediate
from the definition of the group Hilbert algebra and its associated
trace. The trace criterion allows to discuss admissibility in terms
of the central decomposition of the regular representation. In particular,
we present a new proof of the admissibility criteria derived for 
the type I case. In addition we derive admissibility 
criteria which generalize the Wexler-Raz biorthogonality 
relations characterizing dual windows
for Weyl-Heisenberg frames.
\end{abstract}


\section{Introduction}

Given a representation $(\pi,{\cal H}_\pi)$ of a unimodular, separable
locally compact group $G$, we want to discuss the existence and
characterization of vectors giving rise to coherent state expansions on 
${\cal H}_\pi$. For this purpose, a vector $\eta
\in {\cal H}_\pi$ is called {\bf bounded} if the coefficient operator
\[ V_\eta : {\cal H}_\pi \to {\rm L}^2(G)~~,~~ (V_\eta \varphi)(x)
 = \langle \varphi, \pi(x) \eta \rangle ~~\]
is a bounded map.
We are interested in inverting this operator, hence the following notion
is natural: A pair
of bounded vectors $(\eta,\psi)$ is called {\bf admissible} if $V_\eta^* V_\phi$ is
the identity operator on ${\cal H}_\pi$. Note that this property gives rise to
the weak-sense inversion formula
\[
 z = \int_G \langle z, \pi(x) \eta \rangle ~ \pi(x) \psi ~d\mu_G(x)~~, 
\] 
which can be read as a continuous expansion of $z$ in terms of the orbit $\pi(G) \psi \subset
{\cal H}_\pi$. Identities of this type are known as {\bf coherent state expansion} in mathematical physics. 
A single vector $\eta$ is called admissible
if $(\eta,\eta)$ is an admissible pair. It is obvious from the definition that $(\eta,\psi)$
is admissible iff $(\psi,\eta)$ is. In such a case $\eta$ is called the 
{\bf dual vector} of $\psi$.

Admissible vectors were first discussed almost exclusively in connection
with irreducible, so-called {\em discrete series} or {\em square-integrable}
representations
\cite{GrMoPa}. The existence of admissible vectors for these representations
is a fairly straightforward consequence of Schur's Lemma, and slightly more
complicated in the nonunimodular case. Recently exhaustive criteria for the
existence and characterization of admissible vectors were established for the case
that the regular
representation $\lambda_G$ of $G$ is type I, using the Plancherel formula
of the group. However, if we want to include discrete groups in this general
discussion, the type I restriction is rather too rigid. Indeed, a discrete group $G$ has
a type I regular representation iff $G$ itself is type I \cite{Ka}, and the latter
is only the case if $G$ is a finite extension of an abelian normal subgroup \cite{Th}.

Let us now give a short survey of the paper.
It initiated from the idea to replace the decomposition
into irreducibles by the central decomposition of $\lambda_G$, and to
try to come up with criteria using the latter. This however requires an
understanding of how admissible vectors are recognised in terms of their
image under the central decomposition. The main result of this paper, Theorem
\ref{thm:main}, provides the key to this problem, by relating admissibility
to the natural trace on the von Neumann algebra $VN_r(G)$. Since the trace decomposes, 
the problem of characterizing admissible pairs can be translated to the
characterization of tracial pairs for the fibre von Neumann algebras of the
central decomposition (Proposition \ref{prop:adm_fibres}).
As an application of this result we obtain a characterization of admissible
pairs in the case that $\lambda_G$ is type I (Theorem \ref{thm:adm_fibres_typeI}). 
A slightly weaker version of this result had been proved, by somewhat different arguments,
in an earlier paper (\cite[Theorem 1.6]{Fu}).
In the final section we sketch how the trace criterion gives rise
to admissibility criteria in terms of certain orthogonality relations. As a special
case we obtain the Wexler-Raz biorthogonality relations in Gabor analysis.

\section{The group Hilbert algebra, traces and admissible pairs}

Throughout the paper, $G$ denotes a separable unimodular locally compact group, 
and $(\pi,{\cal H}_\pi)$ a (unitary,
strongly continuous) representation of $G$. $\lambda_G$ is the {\bf left regular
representation}, acting on ${\rm L}^2(G)$ by $(\lambda_G(x) f) (y) = f(x^{-1}y)$.
We denote the commuting algebra of $\lambda_G$
as $VN_r(G)$, the {\bf right group von Neumann algebra}.

Our definitions and notations regarding the group Hilbert algebra are taken
from \cite{DiC}.
In order to define the group Hilbert algebra, we let for $f,g \in {\rm L}^2(G)$,
and additionally either $f$ or $g$ in $C_c(G)$,
\[ {\cal U}_f (g) = g \ast f ~~.\]
The (full) group Hilbert algebra consists of all $f \in {\rm L}^2(G)$ for which 
${\cal U}_f$ extends to a bounded operator. Note that these operators then
lie in $VN_r(G)$. Writing $f^*(x) = \overline{f(x^{-1})}$, we note that $V_f = {\cal U}_{f*}
= {\cal U}_f^*$. Hence the bounded vectors are precisely the elements of the full
Hilbert algebra. We note in passing that the full Hilbert algebra contains the
dense subspace ${\rm L}^1(G) \cap {\rm L}^2(G)$ of ${\rm L}^2(G)$.

Let us recall the definition of a trace. Given a von Neumann algebra
${\cal A}$, we let ${\cal A}^+$ denote the cone of positive elements. A mapping
$tr: {\cal A}^+ \to\RR^+ \cup \{ \infty \}$ is called a {\bf trace} if it satisfies the
following two properties:
\begin{itemize}
 \item $tr(S+\alpha T) = tr(S) + \alpha tr(T)$, for all $S,T \in {\cal A}^+$, $\alpha \in \RR^+$.
 Note the conventions $\alpha \infty = \infty$ for $\alpha>0$ and $0 \infty = 0$.
 \item $tr(UTU^*) = tr(T)$, for all $T \in {\cal A}^+$ and all unitary $U \in {\cal A}$.
\end{itemize}
Further relevant properties, that traces may or may not have, are faithfulness, normality 
and semifiniteness; we refer the reader to \cite{DiW} for the definitions. 
$tr$ is {\bf finite} if  $tr({\rm Id}_{\cal H}) = 1$. 
A trace $tr$ uniquely extends to a linear functional on the two-sided ideal
\[ \mathfrak{M}_{tr} = \{ S \in {\cal A} : tr(|S|)< \infty \} ~~,\]
for finite traces this is obviously ${\cal A}$ itself.
We will denote the extension by $tr$ as well. If the trace is normal, 
the associated linear functional is ultra-weakly continuous \cite[III.6, Proposition 1]{DiW}.

The group Hilbert algebra induces a faithful, normal and semifinite trace on $VN_r(G)^+$,
by letting for $T \in VN_r(G)^+$
\[ tr(T) = \left\{ \begin{array}{cc} \| f \|^2 & T = {\cal U}_f^* {\cal U}_f
 \mbox{ for a bounded vector } f \\
 \infty & \mbox{otherwise} \end{array} \right.
\]
We note that for bounded vectors $f,g \in {\rm L}^2(G)$, $V_f^* V_g$ is
in $\mathfrak{M}_{tr}$, with 
\begin{equation}
 \label{eqn:trace_VN} 
tr(V_f^* V_g) = \langle f, g \rangle
\end{equation}

For our arguments it will sometimes be convenient to assume
that ${\cal H}_\pi = {\cal H} \subset {\rm L}^2(G)$ is a closed leftinvariant subspace, on 
which $\pi$ acts by left translation. This is not a restriction, thanks
to the following lemma which collects a few facts about admissible pairs.
We expect most of these statements to be widely known.

\begin{lemma} \label{lem:basics}
\begin{enumerate}
\item[(a)]
 For any bounded vector $\eta$, $V_\eta$ intertwines $\pi$ with $\lambda_G$.
\item[(b)]
 If $\pi$ has an admissible pair $(\eta,\psi)$, then both $V_\eta$ and
 $V_\psi$ are topological embeddings into ${\rm L}^2(G)$. Conversely,
 given a bounded vector $\eta$ such that $V_\eta$ is a topological embedding, 
 a dual vector for $\eta$ is given by $\psi = (V_\eta^* V_\eta)^{-1} \eta$.
 $\psi$ is the unique dual vector with minimal norm.
\item[(c)] There
 exists an admissible pair $(\eta,\psi)$ iff there exists an admissible
 vector.
\item[(d)] If there exists an admissible pair, $\pi$ is
 unitarily equivalent to a subrepresentation of $\lambda_G$.
\end{enumerate}
\end{lemma}
\begin{prf}
 Part (a) is immediate, and (d) then follows from (c). 
 The first statement of (b) is obvious. For the existence of a dual vector, we observe
 that $S=V_\eta^* V_\eta$ is a strictly positive operator with bounded inverse.
 Hence $V_{S^{-1} \eta} = V_\eta \circ S^{-1}$ is bounded,
 i.e., $S^{-1} \psi$ is a bounded vector. The computation 
\[ V_\psi^*V_\eta = V_{S^{-1} \eta}^* V_\eta = S^{-1} V_\eta^* V_\eta =
 {\rm Id}_{{\cal H}_\pi} ~~,\]
 shows that $\psi$ is a dual vector. A similar calculation shows $ V_{S^{-1/2} \eta}^* V_{S^{-1/2} \eta}
 =  {\rm Id}_{{\cal H}_\pi}$, 
 i.e., $(c)$.

 Hence, for the proof of minimality of $\| \psi \|$ (which is the only thing
 left to show), we may assume that ${\cal H}_\pi = {\cal H} \subset {\rm L}^2(G)$,
 and that $\pi$ is left translation on ${\cal H}$. Then
 the commuting algebra $\pi(G)'$ is readily identified as the {\bf
 reduced von Neumann algebra} $\{ p T p : T \in VN_r(G) \} \subset VN_r(G)$, 
 where $p$ denotes the projection onto ${\cal H}$.

 The set of all dual windows is an affine subspace, since the difference of two dual windows 
 is in the linear subspace
 \[ W = \{ x \in {\cal H}_\pi \mbox{ bounded vector } : V_x^* V_\eta = 0 \} ~~.\]
 Hence a dual window of minimal norm is necessarily unique. Now (\ref{eqn:trace_VN}) entails for $x \in W$
 \[ \langle x, \psi \rangle = tr(V_x^* V_\psi) = tr(V_x^* V_{S^{-1} \eta}) = tr(V_x^* V_\eta S^{-1}) 
 = tr(0) = 0 ~~,\]
 hence $\psi \bot W$, and $\| \psi \|$ is indeed minimal.
\end{prf}

The characterization of admissible vectors in terms of the trace
requires one more piece of notation: Given a
particular trace $tr$ on a von Neumann algebra ${\cal A}$, we call a
pair of elements $(\eta,\psi)$ of the underlying Hilbert space {\bf tracial}
if
\[ \forall T \in {\cal A}^+ : tr(T) = \langle T \eta, \psi \rangle ~~.\]

\begin{thm}
\label{thm:main}
 Let ${\cal H} \subset {\rm L}^2(G)$ be a closed, leftinvariant
subspace, with associated leftinvariant projection $p$, and let $\pi$ denote
the restriction of $\lambda_G$ to ${\cal H}$.
\begin{enumerate}
\item[(a)] There exists an admissible pair for ${\cal H}$ iff $tr(p)<\infty$.
\item[(b)] For all pairs $(\eta, \psi) \in {\cal H} \times {\cal H}$ of bounded
 vectors: $(\eta,\psi)$ is admissible iff $(\eta,\psi)$ is tracial for
 $\pi(G)'$.
\end{enumerate}
\end{thm}
\begin{prf}
 For part $(a)$, first assume that there exists an admissible vector $\eta$.
 Hence $p = V_\eta^* V_\eta  = {\cal U}_{\eta^*}^*
 {\cal U}_{\eta^*}$, and thus $tr(p) = \| \eta \|^2 < \infty$.

 Conversely, if $tr(p) < \infty$, then $p =  {\cal U}_{\eta^*}^* {\cal U}_{\eta^*}$,
 for some bounded vector $\eta \in {\rm L}^2(G)$. Now the computation
\[ V_{p\eta}^* V_{p \eta} = p V_\eta^* V_\eta p = p \]
 shows that $p \eta \in {\cal H}$ is admissible.

 For $(b)$, we compute, for any $T = V_g^* V_g$ with $g$ bounded, and for
 any pair $(\eta, \psi)$ of bounded vectors,
 \begin{eqnarray*}
 \langle T \eta, \psi \rangle & = & \langle \eta \ast g^* \ast g, \psi \rangle \\
 & = & \langle \eta \ast g^*, \psi \ast g^* \rangle \\
 & = & \langle g \ast \eta^*, g \ast \psi^* \rangle \\
 & = & \langle g, g \ast \psi^* \ast \eta \rangle \\
 & = & \langle g, (V_\eta^* V_\psi) g \rangle~~.
 \end{eqnarray*}
Hence, assuming that $(\eta,\phi)$ are admissible, we obtain
\[ \langle T \eta, \psi \rangle = \| g \|^2 = tr(T) ~~,\]
as desired. Conversely, assuming traciality of $(\eta,\psi)$, the above
calculation yields 
\[ \| g \|^2 = tr(T) = \langle T \eta, \psi \rangle = \langle g, V_\eta^* V_\psi g \rangle ~~,\]
for all bounded vectors.
By polarization this leads to
\[ \langle h, g \rangle = \langle h, V_\eta^* V_\psi g \rangle ~~,\]
 for all bounded vectors $h,g$, and since these are dense, $V_\eta^* V_\psi = {\rm Id}_{\cal H}$ follows.
\end{prf}

\begin{rem}
 The equivalent conditions from part $(a)$ imply in particular that $\pi(G)'$
 is a finite von Neumann algebra. However, finiteness of $\pi(G)'$ is not
 sufficient, as the following equivalences show:
 \begin{eqnarray}
 \label{eqn:SIN} VN_r(G) \mbox{ is finite } & \Leftrightarrow & G \mbox{ is an SIN-group } \\
 \label{eqn:Discrete} tr({\rm Id}_{{\rm L}^2(G)}) < \infty & \Leftrightarrow & G \mbox{ is discrete }
 \end{eqnarray}
 Here (\ref{eqn:SIN}) is \cite[13.10.5]{DiC}, whereas (\ref{eqn:Discrete}) follows combining
 Theorem \ref{thm:main} (a) with \cite[Proposition 0.4]{Fu}. Recall that SIN-groups
 are defined by having a conjugation-invariant neighborhood-base at unity. Clearly
 this class comprises the locally compact abelian groups, hence for any nondiscrete
 LCA group $VN_r(G)$ is finite, but $tr({\rm Id}_{{\rm L}^2(G)})=\infty$.
\end{rem}

\section{Application to the central decomposition}

In this section we consider the central decomposition of the
regular representation and its use for the characterization of admissible
pairs. In particular we recover the characterization obtained in \cite{Fu}
for the case that $\lambda_G$ is type I. The following facts concerning
the central decomposition of $\lambda_G$ can be found in \cite[18.7.7,18.7.8]{DiC}.

Let $\check{G}$ denote the space of quasi-equivalence classes of factor representations
of $G$, endowed with the natural Borel structure. Then there exists a standard positive
measure $\nu_G$ on $\check{G}$, and a measurable field of factor representations $\rho_\sigma
\in \sigma$, for $\nu_G$-almost every $\sigma$, such that
\[
 \lambda_G \simeq \int_{\check{G}}^\oplus \rho_\sigma d\nu_G(\sigma)~~.
\]
The operator effecting the unitary equivalence is called the {\bf Plancherel transform}.
Moreover, the direct integral provides a decomposition of $VN_r(G)$ and the natural trace:
\[ VN_r(G) =  \int_{\check{G}}^\oplus {\cal A}_\sigma d\nu_G(\sigma) \]
for a measurable family of von Neumann algebras ${\cal A}_\sigma$ on ${\cal H}_\sigma$, as well as
\[ tr(T) = \int_{\check{G}} tr_\sigma(T_\sigma) d\nu_G(\sigma)~~,
\]
when $(T_\sigma)_{\sigma \in \check{G}}$ denotes the operator field corresponding
to $T$ under the central decomposition, and $tr_\sigma$ is a faithful normal, semifinite
trace on the factor ${\cal A}_\sigma$, which exists for $\nu_G$-almost every $\sigma$.
In particular, $\nu_G$-almost every ${\cal A}_\sigma$ is of type I or II.

Now admissibility is easily translated to traciality in the fibres. In the following,
$(\widehat{\psi}_{\sigma})_{\sigma \in \check{G}}$ denotes the Plancherel
transform of $\psi \in {\rm L}^2(G)$.
\begin{prop}
\label{prop:adm_fibres}
 Let $\pi$ denote the restriction of $\lambda_G$ to a closed, leftinvariant
 subspace ${\cal H} \subset {\rm L}^2(G)$. Let $P$ denote the projection onto ${\cal H}$,
 then $P$ decomposes into a measurable field of projections $\widehat{P}_\sigma$,
 and $\pi(G)'$ decomposes under
 the central decomposition into the von Neumann algebras ${\cal C}_\sigma = \widehat{P}_\sigma
 {\cal A}_\sigma \widehat{P}_\sigma$.
\begin{enumerate}
\item[(a)]  For bounded $\eta,\psi \in {\cal H}$, we have
 \[ (\eta,\psi) \mbox{ is admissible for ${\cal H}$ } \Leftrightarrow  
  (\widehat{\eta}_\sigma, \widehat{\psi}_\sigma) \mbox{ is tracial for } {\cal C}_\sigma ~~(\nu_G \mbox{a.e.}) \]
\item[(b)] ${\cal H}$ has an admissible pair of vectors iff $\int_{\check{G}} tr(\widehat{P}_\sigma) d\nu_G(\sigma)
< \infty$. In particular, almost all $C_\sigma$ are finite von Neumann algebras.
\end{enumerate}
\end{prop}
The (potential) use of the proposition consists in the fact that we only need to characterize
tracial pairs for factor representations. Unfortunately, we are not aware of any explicit criteria
for tracial vectors associated to type II factors. For type I factors, though, they are easily derived,
as the proof of the next theorem shows. Note that if $\lambda_G$ is type I, the fibre spaces
in the central decomposition are just the Hilbert-Schmidt spaces ${\cal B}_2({\cal H}_\sigma)$,
where $\sigma$ runs through the unitary dual, and $\lambda_G$ decomposes into left action on
${\cal B}_2({\cal H}_\sigma)$ via $\sigma$ \cite[18.8]{DiC}. 
\begin{thm} \label{thm:adm_fibres_typeI}
 Let $G$ be unimodular with $\lambda_G$ type I. Let ${\cal H} \subset {\rm L}^2(G)$ be
 a leftinvariant subspace. Then there exists a measurable field of projections $\widehat{P}_\sigma$
 on ${\cal H}_\sigma$ such that
  \[ P \simeq \int_{\widehat{G}}^\oplus 1 \otimes \widehat{P}_\sigma ~ d\nu_G(\sigma) ~~.\]
 \begin{enumerate}
 \item[(a)] $(\eta,\psi)$ is admissible $\Leftrightarrow$ for $\nu_G$-almost every $\sigma \in \widehat{G}$ :
 $\widehat{\psi}_\sigma^* \widehat{\eta}_\sigma = \widehat{P}_\sigma$.
 \item[(b)] There exists an admissible vector for ${\cal H}$ iff
 \[
  \nu_{\cal H} = \int_{\widehat{G}} {\rm rank}(\widehat{P}_\sigma) d\nu_G(\sigma) < \infty.
 \]
 \end{enumerate}
\end{thm}
\begin{prf}
The existence of the $\widehat{P}_\sigma$ follows from the type I property.
In the following it is convenient to use tensor-product notation for rank-one operators,
i.e., $x \otimes y$ denotes the operator $z \mapsto \langle z,y \rangle x$.
Given a fixed $\sigma \in \widehat{G}$,
the elements of ${\cal K} = {\cal B}_2({\cal H}_\sigma) \circ \widehat{P}_\sigma$ can be
written uniquely as $\eta = \sum_{i \in I} \eta_i \otimes e_i$, where $(e_i)_{i \in I}$ is a fixed
orthogonal basis of $\widehat{P}({\cal H}_\sigma)$. It follows that ${\cal K}$
is conveniently identified with ${\cal H}_\sigma \otimes \ell^2(I)$. In this
identification the left action of $\sigma$ on ${\cal K}$ becomes $\sigma \otimes 1$.
Moreover, the commuting algebra is easily identified with $1 \otimes {\cal B}_2(\ell^2(I))$,
and its trace is the usual operator trace.

A weak-operator dense subspace of ${\cal B}(\ell^2(I))$ is spanned by the operators
$e_{i,k} = \delta_i \otimes \delta_k$, where $\delta_i \in \ell^2(I)$
denotes the usual Kronecker-$\delta$ concentrated at $i$.
Now, given $\widehat{\eta}_\sigma=\sum_{i \in I} \eta_i \otimes e_i$ and 
$\widehat{\psi}_\sigma = \sum_{i \in I} \psi_i \otimes e_i$, we compute
\[ tr(e_{i,k}) = \delta_{i,k} \]
and
\[ \langle (1 \otimes e_{i,k}) \eta, \psi \rangle = \langle \eta_i, \psi_k \rangle ~~,\] 
whence we obtain the following
traciality condition
\begin{eqnarray*}
 (\widehat{\eta}_\sigma,\widehat{\psi}_\sigma) \mbox{ tracial } &
 \Longleftrightarrow & \forall i,k ~:~ \langle \eta_i, \psi_k \rangle = \delta_{i,k} \\
 & \Longleftrightarrow & \left(\sum_{i \in I} \psi_i \otimes e_i\right)^*
 \left( \sum_{i \in I} \eta_i \otimes e_i \right) = \widehat{P}_\sigma ~~,
\end{eqnarray*}
which proves part $(a)$. Part $(b)$ follows easily from $(a)$, see \cite{Fu}.
\end{prf}

We wish to point out that the admissibility criteria, however abstract they may appear,
have been made explicit for certain classes of representations, in particular for
multiplicity-free representations. See \cite{FuMa} for a discussion of quasiregular
representations of semidirect product groups, and \cite{Fu_Ch} for a treatment of
Weyl-Heisenberg frames with integer sampling ratio.

Another interesting class of representations are the factor subrepresentations
of the regular representation, i.e., the atoms in the central decomposition, and
the elements of their quasi-equivalence classes.
These representations were already considered in \cite{Ro}, though not with a view
to constructing admissible vectors.
\begin{cor} Let $\pi$ be a factor representation.
\begin{enumerate}
\item[(a)] $\pi$ has admissible vectors iff $\pi$ is equivalent to a subrepresentation
of $\lambda_G$, and $\pi(G)'$ is a finite von Neumann algebra. In particular, $\pi$ has
either type I or II, and there exists a faithful, finite, normal trace $tr$ on $\pi(G)'$, unique
up to normalization.
\item[(b)] The trace on $\pi(G)'$ can be normalized in such a way that the following equivalence
holds:
\[ (\eta,\psi) \mbox{ is admissible } \Leftrightarrow (\eta,\psi) \mbox{ is tracial } ~~.\]
\end{enumerate}
\end{cor}

\section{Checking admissibility using biorthogonality relations}

While the discussion of the type I case shows that the characterization of admissible vectors
via the trace on the commutant can be used to some effect, in the general case the merits are
much less obvious. In this section we sketch a procedure to arrive at more concrete
necessary and sufficient conditions for admissible pairs, in terms of certain scalar products.
We will then demonstrate that the Wexler-Raz biorthogonality relations are a special instance
of this approach. For the formulation of the admissibility conditions, we require
\begin{itemize}
 \item A family $(T_i)_{i \in I} \subset \pi(G)'$ spanning a weak-operator dense subspace of
 $\pi(G)'$. Recall that the density requirement means that for each $S \in \pi(G)'$ there exists
 a net $(S_j)_{j \in J}$ in the span such that for all pairs $y,z \in {\cal H}_\pi$ we have
 $\langle S_j y, z \rangle \to \langle S y, z \rangle$.
 \item An admissible pair $(\eta_0, \psi_0)$.
\end{itemize}
Then for a pair of bounded vectors $(\eta,\psi)$ we have the following equivalence:
\begin{equation} \label{eqn:gen_WR}
(\eta, \psi) \mbox{ is admissible } \Longleftrightarrow \forall i \in I ~:~ \langle T_i \eta, \psi \rangle
 = \langle T_i \eta_0, \psi_0 \rangle~~.
\end{equation}
The proof of the condition is immediate from the assumptions and Theorem \ref{thm:main}.
The criterion is explicit as soon as the $T_i$ and the admissible pair $(\eta_0,\psi_0)$
are known explicitly. Clearly, generators are preferable which provide particularly
simple relations.

Let us now derive admissibility criteria in the context of Weyl-Heisenberg frames. These frames 
are obtained by picking a window function
$\eta \in {\rm L}^2(\RR)$ and translating it along a time-frequency lattice
$\Gamma$. The shifts are described in terms of the operators
\[ T_x : f \mapsto  f(\cdot - x) ~~,~~ M_\omega : f \mapsto e^{2 \pi i \omega \cdot}
 f ~~.\]
For the following we fix $\alpha, \beta >0$, with $\alpha \beta \le 1$. This
is a well-known necessary and sufficient condition for the existence of Weyl-Heisenberg frames.
The sufficiency is proved by the admissible vector $\eta_0$ given below, for necessity confer,
among others, \cite{Ba,Be,Ri}. 
Given $\eta \in {\rm L}^2(\RR)$, we wish to decide whether the family
\[ \{ M_{\alpha m} T_{\beta n} \eta : m,n \in \ZZ \} \]
constitutes a frame of ${\rm L}^2(\RR)$. Recall that the latter
property means that the coefficient map
\[ {\cal T}_{f;\alpha,\beta}: f \mapsto \left( \langle f, M_{\alpha m} T_{\beta n} \eta \rangle
 \right)_{m,n \in \ZZ} \]
defines a topological embedding ${\rm L}^2(\RR) \hookrightarrow \ell^2(\ZZ^2)$.
As in the proof of Lemma \ref{lem:basics} we see that $f$ generates a Weyl-Heisenberg
frame iff there exists a {\bf dual window} $g$ generating a Weyl-Heisenberg frame
and satisfying ${\cal T}_{g;\alpha,1}^* {\cal T}_{f;\alpha,\beta} = {\rm Id}$.

Since the time-frequency shifts $M_{\alpha m} T_n$ do not constitute a group of
operators, the group-theoretic interpretation of the problem requires a slight
detour in the form of the next lemma. Note that $\TT$ denotes the set
of complex numbers with modulus one.

\begin{lemma}
 Define the group $G = \ZZ \times \ZZ \times \TT$, with group law
 \[ (m,n,z) (m',n',z') = (m+m',n+n',zz' e^{-2 \pi i \alpha \beta m' n}) ~~.\]
 $G$ acts on ${\rm L}^2(\RR)$ via the representation
 \[ \pi(m,n,z) = M_{\alpha m} T_{\beta n} z ~~.\]
 For all $(f,g) \in {\rm L}^2(\RR)$ with ${\cal T}_{f;\alpha,\beta},{\cal T}_{g;\alpha,\beta}$
 bounded, $f$ generates a Weyl-Heisenberg frame with dual window $g$
 iff $(f,g)$ is an admissible pair for $\pi$. 
\end{lemma}
\begin{prf}
 The statements concerning $G$ and $\pi$ are immediate from the definitions. 
 For the last statement, observe that 
 \begin{eqnarray*} V_g^* V_f h & = & \int_\TT \sum_{m,k \in \ZZ} \langle h, \pi(m,k,z) f \rangle
 ~\pi(m,k,z) g ~dz \\
 & = & \int_\TT \sum_{m,k \in \ZZ} \langle h, \pi(m,k,0) f \rangle ~ \pi(m,k,0) g~dz \\
 & = & {\cal T}_{g;\alpha,\beta}^* {\cal T}_{f;\alpha,\beta} ~~. \end{eqnarray*}
\end{prf}

\begin{rem}
 The representation $\pi$ is type I iff $\alpha \beta$ is rational. 
 For the only-if part confer \cite[Remark 2.]{Ba}, whereas the if-part
 follows from the fact that the group itself is type I if $\alpha \beta$ is rational
 (a straightforward application of Mackey's theory).  
 In the case where $1/(\alpha \beta) \in \ZZ$, there exist admissibility criteria
 which employ the so-called Zak transform; here the representation is even multiplicity-free.
 See \cite{Fu_Ch} for an interpretation of the Zak transform criterion in the light of
 Theorem \ref{thm:adm_fibres_typeI}
\end{rem}

Following the general procedure sketched above, we now observe that
\begin{itemize}
\item $\eta_0 = \sqrt{\alpha} \chi_{[0,\beta)}$ is an admissible vector \cite{DaGr}.
\item The {\bf commuting lattice}
 \[ \Lambda_c = \{ M_{m/\beta} T_{n/\alpha} : m, n \in \ZZ \} 
 \] 
 generates a weak-operator dense subspace of $\pi(G)'$ (\cite[Appendix 6.1]{DaLaLa}).
\end{itemize}
Hence, after verifying that 
\[ \langle M_{m/\beta} T_{n/\alpha} \eta_0, \eta_0 \rangle = \alpha \beta
 \delta_{m,0} \delta_{n,0}~~, \]
we obtain the {\bf Wexler-Raz biorthogonality relations} as a special case of
 (\ref{eqn:gen_WR}):
\begin{cor}
 Let $g,\gamma$ be such that ${\cal T}_{g; \alpha,\beta}, {\cal T}_{\gamma; \alpha,\beta}$ are bounded.
 Then $\gamma$ is a dual window for $g$ iff
 \begin{equation} \label{eqn:WR_bior}
 \langle  M_{m/\beta} T_{n/\alpha} \gamma, g \rangle = \alpha \beta
 \delta_{m,0} \delta_{n,0} ~~.\end{equation} 
\end{cor}

\begin{rem}
 A more general ``Wexler-Raz-relation'' is  
 \begin{equation} \label{eqn:WR}
 {\cal T}_{f;\alpha, \beta}^* {\cal T}_{g;\alpha, \beta} h = \frac{1}{\alpha \beta} {\cal T}_{h;1/\beta,1/\alpha}^*
 {\cal T}_{g;1/\beta,1/\alpha} f 
 \end{equation}
 proved for suitable $f,g,h$ in \cite{DaLaLa}. (\ref{eqn:WR}) is easily seen to imply
 (\ref{eqn:WR_bior}). It is not clear whether (\ref{eqn:WR}) has
 a counterpart in the general setting.  
\end{rem}

\section*{Concluding remarks}

Von Neumann algebra techniques have been used previously for establishing
criteria for the existence of cyclic and/or admissible vectors, see for instance
\cite{Ba,Ri}. In particular the coupling constant has proved to be a powerful
tool for existence results, see \cite{Be,Ri}. However, these techniques seem to
be of limited use for the explicit construction of admissible vectors. By contrast, 
this paper aims at providing criteria for these vectors,
though it is clear that much remains to be done to make these criteria work. The authors of 
\cite{DaLaLa} used the trace on the commuting algebra in the Weyl-Heisenberg
frame context, but did not point out the close connection to admissibility.


\section*{References}
\begin{thebibliography}{99}
  \bibitem{Ba}{L. Baggett, {\em Processing a radar signal and representations of the
 discrete Heisenberg group}, Coll. Math. {\bf 60/61} (1990), 195-203.}
  \bibitem{Be}{M.B. Bekka, {\em Square integrable representations, von Neumann algebras
 and an application to Gabor analysis}, Preprint, 2002.}
  \bibitem{DaGr}{I. Daubechies and A. Grossmann, {\em Painless nonorthogonal expansions,}
 J. Math. Phys. {\bf 27} (1986), 1271-1283.}
  \bibitem{DaLaLa}{I. Daubechies, H.J. Landau and Z. Landau, {\em Gabor time-frequency
  lattices and the Wexler-Raz identity,} J. Fourier Analysis and Applications {\bf 1} (1995),
  437-478.}
  \bibitem{DiC}{J. Dixmier, {\em $C^*$-Algebras.} North Holland, Amsterdam, 1977.}
  \bibitem{DiW}{J. Dixmier, {\em Von Neumann Algebras.} North Holland,
 Amsterdam, 1981.}
 \bibitem{Fu}{H. F\"uhr, {\em Admissible vectors for the regular representation,}
 Proc. AMS {\bf 130} (2002), 2959-2970.}
 \bibitem{Fu_Ch}{H. F\"uhr, {\em Plancherel transform criteria for Weyl-Heisenberg frames with 
 integer oversampling}, submitted, available as math.FA/0206309}
 \bibitem{FuMa}{H. F\"uhr and M. Mayer, {\em Continuous wavelet transform from semidirect products:
 Cyclic representations and Plancherel measure}, J. Fourier Analysis and Applications, to appear.}
 \bibitem{GrMoPa}{A. Grossmann, J. Morlet and T. Paul, {\em Transforms associated to square integrable
 group representations I: General results}, J. Math. Phys. {\bf 26} (1985), 2473-2479.} 
 \bibitem{Ka}{E. Kaniuth, {\em Der Typ der regul\"aren Darstellung diskreter Gruppen,} Math. Ann.
 {\bf 182} (1969), 334-339.}
 \bibitem{Ri}{M. Rieffel, {\em Von Neumann algebras associated to pairs of lattices in Lie groups,}
 Math. Ann. {\bf 257} (1981), 403-418.}
 \bibitem{Ro}{J. Rosenberg, {\em Square-integrable factor representations of
 locally compact groups,} Trans. AMS {\bf 261} (1978), pp. 1-33.}
 \bibitem{Th}{E. Thoma, {\em Eine Charakterisierung diskreter Gruppen vom Typ I}, 
 Invent. Math. {\bf 6} (1968), 190-196.}
\end{thebibliography}
\end{document}